\documentclass[5p,number]{elsarticle}
\usepackage{extarrows,paralist}
\usepackage{mhchem} 
\usepackage{siunitx} 
\usepackage{graphicx} 
\usepackage{amsmath,graphicx}
\usepackage{mathtools}		
\usepackage{mathtools}
\usepackage{graphicx}
\usepackage{float}
\usepackage{caption}
\usepackage{subcaption}
\usepackage{color}
\usepackage{placeins}
\usepackage{float}
\usepackage{tabularx,colortbl}
\usepackage{amsmath}
\usepackage{amssymb}
\newcommand*{\QEDB}{\hfill\ensuremath{\square}}%
\usepackage{listings}
\usepackage{epstopdf}
\usepackage{amssymb}
\usepackage{algorithm}
\usepackage{algpseudocode}

\usepackage{enumitem}
\allowdisplaybreaks
\newtheorem{thm}{Theorem}
\newtheorem{lem}{Lemma}
\newtheorem{pf}{Proof}
\newtheorem{assum}{Assumption}
\newtheorem{prop}{Proposition}

\def\mbb{\mathbb}
\def\mb{\mathbf}
\def\mc{\mathcal}

\usepackage{amssymb}

\begin{document}
	\begin{frontmatter}
		\title{On the Distributed Optimization over Directed Networks\tnoteref{t1}}
		\tnotetext[t1]{This work has been partially supported by an NSF Career Award \# CCF-1350264.}
		
		\author{Chenguang Xi}\ead{chenguang.xi@tufts.edu}   
		\author{Qiong Wu}\ead{qiong.wu@tufts.edu}
		\author{Usman A. Khan\corref{cor1}}\ead{khan@ece.tufts.edu}               
		\cortext[cor1]{Corresponding author}
		
		\address{Department of Electrical and Computer Engineering, Tufts University, 161 College Ave., Medford, MA, 02155, USA}  
		\address{Department of Mathematics, Tufts University, 503 Boston Ave., Medford, MA 02155, USA}  
		\begin{abstract}
		In this paper, we propose a distributed algorithm, called Directed-Distributed Gradient Descent (D-DGD), to solve multi-agent optimization problems over \emph{directed} graphs. Existing algorithms mostly deal with similar problems under the assumption of undirected networks, i.e., requiring the weight matrices to be doubly-stochastic. The row-stochasticity of the weight matrix guarantees that all agents reach consensus, while the column-stochasticity ensures that each agent's local gradient contributes equally to the global objective. In a directed graph, however, it may not be possible to construct a doubly-stochastic weight matrix in a distributed manner. We overcome this difficulty by augmenting an additional variable for each agent to record the change in the state evolution. In each iteration, the algorithm simultaneously constructs a row-stochastic matrix and a column-stochastic matrix instead of only a doubly-stochastic matrix. The convergence of the new weight matrix, depending on the row-stochastic and column-stochastic matrices, ensures agents to reach both consensus and optimality. The analysis shows that the proposed algorithm converges at a rate of~$O(\frac{\ln k}{\sqrt{k}})$, where~$k$ is the number of iterations.
		\end{abstract}
		
		\begin{keyword}
			
			
			Distributed optimization; multi-agent networks; directed graphs; distributed gradient descent.             
		\end{keyword}
		
	\end{frontmatter}

\section{Introduction}\label{s1}
Distributed computation and optimization,~\cite{distributed_Tsitsiklis,distributed_Tsitsiklis2}, has received significant recent interest in many areas, e.g., multi-agent networks,~\cite{distributed_Boyd}, model predictive control,~\cite{distributed_Necoara,distributed_Necoara2}, cognitive networks,~\cite{distributed_Mateos}, source localization,~\cite{distributed_Rabbit}, resource scheduling,~\cite{distributed_Chunlin}, and message routing,~\cite{distributed_Neglia}. The related problems can be posed as the minimization of a sum of objectives,~$\sum_{i=1}^{n}f_i(\mb{x})$, where~$f_i:\mbb{R}^p\rightarrow\mbb{R}$ is a private objective function at the $i$th agent. There are two general types of distributed algorithms to solve this problem. The first type is a gradient based method, where at each iteration a gradient related step is \mbox{calculated}, followed by averaging with neighbors in the network, e.g., the Distributed Gradient Descent (DGD),~\cite{uc_Nedic}, and the distributed dual-averaging method,~\cite{cc_Duchi}. The main advantage of these methods is computational simplicity. The second type of distributed algorithms are based on augmented Lagrangians, where at each iteration the primal variables are solved to minimize a Lagrangian related function, followed by updating the dual variables accordingly, e.g., the Distributed Alternating Direction Method of Multipliers (D-ADMM),~\cite{ADMM_Mota,ADMM_Wei, ADMM_Shi}. The latter type is preferred when agents can solve the local optimization problem efficiently. All proposed distributed algorithms,~\cite{uc_Nedic,cc_Duchi,ADMM_Mota,ADMM_Wei,ADMM_Shi}, assume undirected graphs. The primary reason behind assuming the undirected graphs is to obtain a doubly-stochastic weight matrix. The row-stochasticity of the weight matrix guarantees that all agents reach consensus, while the column-stochasticity ensures optimality, i.e., each agent's local gradient contributes equally to the global objective.

In this paper, we propose a gradient based method solving distributed optimization problem over~\emph{directed} graphs, which we refer to as the Directed-Distributed Gradient Descent (D-DGD). Clearly, a directed topology has broader applications in contrast to undirected graphs and may further result in reduced communication cost and simplified topology design. We start by explaining the necessity of weight matrices being doubly-stochastic in existing gradient based method, e.g., DGD. In the iteration of DGD, agents will not reach consensus if the row sum of the weight matrix is not equal to one. On the other hand, if the column of the weight matrix does not sum to one, each agent will contribute differently to the network. Since doubly-stochastic matrices may not be achievable in a directed graph, the original methods, e.g., DGD, no longer work. We overcome this difficulty in a directed graph by augmenting an additional variable for each agent to record the state updates. In each iteration of the D-DGD algorithm, we simultaneously construct a row-stochastic matrix and a column-stochastic matrix instead of only a doubly-stochastic matrix. We give an intuitive explanation of our proposed algorithm and further provide convergence and convergence rate analysis as well.

In the context of directed graphs, related work has considered distributed gradient based algorithms,~\cite{opdirect_Nedic,opdirect_Nedic2,opdirect_Tsianous,opdirect_Tsianous2,opdirect_Tsianous3}, by combining gradient descent and push-sum consensus. The push-sum algorithm,~\cite{ac_directed0,ac_directed}, is first proposed in consensus problems\footnote{See,~\cite{c_Jadbabaie,c_Reynolds,c_Saber,c_Saber2,c_Saber3,c_Xiao}, for additional information on average consensus problems.} to achieve average-consensus given a column-stochastic matrix. The idea is based on computing the stationary distribution (the left eigenvector of the weight matrix corresponding to eigenvalue 1) for the Markov chain characterized by the multi-agent network and canceling the imbalance by dividing with the left-eigenvector. The algorithms in ~\cite{opdirect_Nedic,opdirect_Nedic2,opdirect_Tsianous,opdirect_Tsianous2,opdirect_Tsianous3} follow a similar spirit of push-sum consensus and propose nonlinear (because of division) methods. In contrast, our algorithm follows linear iterations and does not involve any division.

The remainder of the paper is organized as follows. In Section~\ref{s2}, we provide the problem formulation and show the reason why DGD fails to converge to the optimal solution over directed graphs. We subsequently present the D-DGD algorithm and the necessary assumptions. The convergence analysis of the D-DGD algorithm is studied in Section~\ref{s3}, consisting of agents' consensus analysis and optimality analysis. The convergence rate analysis and numerical experiments are presented in Sections~\ref{s4} and~\ref{s5}. Section~\ref{s6} contains concluding remarks.

\textbf{Notation:} We use lowercase bold letters to denote vectors and uppercase italic letters to denote matrices. We denote by~$[\mb{x}]_i$ the~$i$th component of a vector~$\mb{x}$, and by~$[A]_{ij}$ the $(i,j)$th element of a matrix,~$A$. An~$n$-dimensional vector with all elements equal to one (zero) is represented by~$\mb{1}_n$ ($\mb{0}_n$). The notation~$0_{n\times n}$ represents an~$n\times n$ matrix with all elements equal to zero. The inner product of two vectors~$\mb{x}$ and~$\mb{y}$ is~$\langle\mb{x},\mb{y}\rangle$. We use~$\|\mb{x}\|$ to denote the standard Euclidean norm.

\section{Problem Formulation}\label{s2}
Consider a strongly-connected network of~$n$ agents communicating over a directed graph,~$\mc{G}=(\mc{V},\mc{E})$, where~$\mc{V}$ is the set of agents, and~$\mc{E}$ is the collection of ordered pairs,~$(i,j),i,j\in\mc{V}$, such that agent~$j$ can send information to agent~$i$. Define~$\mc{N}_i^{{\scriptsize \mbox{in}}}$ to be the collection of in-neighbors, i.e., the set of agents that can send information to agent~$i$. Similarly,~$\mc{N}_i^{{\scriptsize \mbox{out}}}$ is defined as the out-neighborhood of agent~$i$, i.e., the set of agents that can receive information from agent~$i$. We allow both~$\mc{N}_i^{{\scriptsize \mbox{in}}}$ and~$\mc{N}_i^{{\scriptsize \mbox{out}}}$ to include the node~$i$ itself. Note that in a directed graph~$\mc{N}_i^{{\scriptsize \mbox{in}}}\neq\mc{N}_i^{{\scriptsize \mbox{out}}}$, in general. We focus on solving a convex optimization problem that is distributed over the above network. In particular, the network of agents cooperatively solve the following optimization problem:
\begin{align}
\mbox{P1}:
\quad\mbox{min  }&f(\mb{x})=\sum_{i=1}^nf_i(\mb{x}),\nonumber
\end{align}
where each~$f_i:\mbb{R}^p\rightarrow\mbb{R}$ is convex, not necessarily differentiable, representing the local objective function at agent~$i$.
\begin{assum}\label{asp}
In order to solve the above problem, we make the following assumptions:
	\begin{enumerate}[label=(\alph*)]
		\item The agent graph,~$\mc{G}$, is strongly-connected.
		\item Each local function,~$f_i:\mbb{R}^p\rightarrow\mbb{R}$, is convex,~$\forall i\in\mc{V}$.
        \item The solution set of Problem P1 and the corresponding optimal value exist. Formally, we have
		\begin{align}
		\mb{x}^*\in\mc{X}^*=\left\{\mb{x}|f(\mb{x})=\min_{\mb{y}\in\mbb{R}^p}f(\mb{y})\right\},f^*=\min f(\mb{x}).\nonumber
		\end{align}
		\item The sub-gradient,~$\nabla f_i(\mb{x})$, is bounded:
		\[
		\|\nabla f_i(\mb{x})\|\leq D,
		\]
		for all~$\mb{x}\in\mbb{R}^p, i\in\mc{V}$.
	\end{enumerate}
\end{assum}
The Assumptions~\ref{asp} are standard in distributed optimization, see related literature,~\cite{cc_nedic}, and references therein. Before describing our algorithm, we first recap the DGD algorithm,~\cite{uc_Nedic}, to solve P1 in an undirected graph. This algorithm requires doubly-stochastic weight matrices. We analyze the influence to the result of the DGD when the weight matrices are \emph{not} doubly-stochastic.	

\subsection{Distributed Gradient Descent}
Consider the Distributed Gradient Descent (DGD),~\cite{uc_Nedic}, to solve P1. Agent~$i$ updates its estimate as follows:
\begin{equation}\label{DGD}
\mb{x}^{k+1}_i=\sum_{j=1}^nw_{ij}\mb{x}^{k}_j-\alpha_{k}\nabla f_i^k,
\end{equation}
where~$w_{ij}$ is a non-negative weight such that~$W=\{w_{ij}\}$ is doubly-stochastic. The scalar,~$\alpha_{k}$, is a diminishing but non-negative step-size, satisfying the persistence conditions,~\cite{kushner2003stochastic,cc_lobel}:~$\sum_{k=0}^\infty\alpha_k=\infty$,~$\sum_{k=0}^\infty\alpha_k^2<\infty$, and the vector~$\nabla f_i^k$ is a sub-gradient of~$f_i$ at~$\mb{x}^k_i$. For the sake of argument, consider~$W$ to be row-stochastic but not column-stochastic. Clearly,~$\mb{1}$ is a right eigenvector of~$W$, and let~$\boldsymbol{\pi}=\{\pi_i\}$ be its left eigenvector corresponding to eigenvalue~$1$.
Summing over~$i$ in Eq.~\eqref{DGD}, we get
\begin{align}\label{dps_row}
\widehat{\mb{x}}^{k+1}&\triangleq\sum_{i=1}^n\pi_i\mb{x}^{k+1}_i,\nonumber\\
&=\sum_{j=1}^n\left(\sum_{i=1}^n\pi_iw_{ij}\right)\mb{x}^{k}_j-\alpha_{k}\sum_{i=1}^n\pi_i\nabla f_i(\mb{x}_i^k),\nonumber\\
&=\widehat{\mb{x}}^{k}-\alpha_{k}\sum_{i=1}^n\pi_i\nabla f^k_i,
\end{align}
where~$\pi_j=\sum_{i=1}^{n}\pi_iw_{ij},\forall i,j$. If we assume that the agents reach an agreement, then Eq.~\eqref{dps_row} can be viewed as an inexact (central) gradient descent (with~$\sum_{i=1}^n\pi_i\nabla f_i(\mb{x}^k_i)$ instead of~$\sum_{i=1}^n\pi_i\nabla f_i(\widehat{\mb{x}}^k)$) minimizing a new objective, $\widehat{f}(\mb{x})\triangleq\sum_{i=1}^n\pi_if_i(\mb{x})$. As a result, the agents reach consensus and converge to the minimizer of~$\widehat{f}(\mb{x})$.

Now consider the weight matrix,~$W$, to be column-stochastic but not row-stochastic. Let~$\overline{\mb{x}}^k$ be the average of agents estimates at time~$k$,  then Eq.~\eqref{DGD} leads to
\begin{align}\label{dps_column}
\overline{\mb{x}}^{k+1}&\triangleq\frac{1}{n}\sum_{i=1}^n\mb{x}^{k+1}_i,\nonumber\\
&=\frac{1}{n}\sum_{j=1}^n\left(\sum_{i=1}^nw_{ij}\right)\mb{x}^{k}_j-\frac{\alpha_{k}}{n}\sum_{i=1}^n\nabla f_i(\mb{x}_i^k),\nonumber\\
&=\overline{\mb{x}}^{k}-\left(\frac{\alpha_{k}}{n}\right)\sum_{i=1}^n\nabla f^k_i.
\end{align}
Eq.~\eqref{dps_column} reveals that the average,~$\overline{\mb{x}}^k$, of agents estimates follows an inexact (central) gradient descent ($\sum_{i=1}^n\nabla f_i(\mb{x}^k_i)$ instead of $\sum_{i=1}^n\nabla f_i(\overline{\mb{x}}^k)$) with stepsize~$\alpha^k/n$, thus reaching the minimizer of~$f(\mb{x})$. Despite the fact that the average,~$\overline{\mb{x}}^k$, reaches the optima,~$\mb{x}^*$, of~$f(\mb{x})$, the optima is not achievable for each agent because consensus can not be reached with a matrix that is not necessary row-stochastic.

Eqs.~\eqref{dps_row} and~\eqref{dps_column} explain the importance of doubly-stochastic matrices in consensus-based optimization. The row-stochasticity guarantees all of the agents to reach a consensus, while column-stochasticity ensures each local gradient to contribute equally to the global objective.

\subsection{Directed-Distributed Gradient Descent (D-DGD)}
From the above discussion, we note that reaching a consensus requires the right eigenvector (corresponding to eigenvalue~$1$) to lie in~$\mbox{span}\{\mb{1}_n\}$, and minimizing the \mbox{global} objective requires the corresponding left eigenvector to lie in~$\mbox{span}\{\mb{1}_n\}$. Both the left and right eigenvectors of a doubly-stochastic matrix are~$\mb{1}_n$, which, in general, is not possible in directed graphs. In this paper, we introduce \emph{Directed-Distributed Gradient Descent} (D-DGD) that overcomes the above issues by augmenting an additional variable at each agent and thus constructing a new weight matrix,~$W\in\mbb{R}^{2n\times 2n}$, whose left and right eigenvectors (corresponding to eigenvalue~$1$) are in the form:~$[\mb{1}_n^\top,\mb{v}^\top]$ and~$[\mb{1}_n^\top,\mb{u}^\top]^\top$. Formally, we describe D-DGD as follows.

At $k$th iteration, each agent,~$j\in\mc{V}$, maintains two vectors:~$\mb{x}_j^k$ and~$\mb{y}_j^k$, both in~$\mbb{R}^p$. Agent~$j$ sends its state estimate,~$\mb{x}_j^k$, as well as a weighted auxiliary variable,~$b_{ij}\mb{y}_j^k$, to each out-neighbor,~$i\in\mc{N}_j^{{\scriptsize \mbox{out}}}$, where~$b_{ij}$'s are such that:
\begin{equation*}
b_{ij}=\left\{
\begin{array}{rl}
>0,&i\in\mc{N}_j^{{\scriptsize \mbox{out}}},\\
0,&\mbox{otw.},
\end{array}
\right.
\qquad
\sum_{i=1}^nb_{ij}=1,\forall j.
\end{equation*}
Agent~$i$ updates the variables,~$\mb{x}_i^{k+1}$ and~$\mb{y}_i^{k+1}$, with the information received from its in-neighbors,~$j\in\mc{N}_i^{{\scriptsize \mbox{in}}}$, as follows:
\begin{subequations}\label{alg1}
\begin{align}
\mb{x}_i^{k+1}&=\sum_{j=1}^na_{ij}\mb{x}_j^k+\epsilon\mb{y}_i^k-\alpha_k\nabla f_i(\mb{x}_i^k),\label{alg1a}\\
\mb{y}_i^{k+1}&=\mb{x}_i^k-\sum_{j=1}^na_{ij}\mb{x}_j^k+\sum_{j=1}^nb_{ij}\mb{y}_j^k-\epsilon\mb{y}_i^k,\label{alg1b}
\end{align}
\end{subequations}
where:
\begin{equation*}
a_{ij}=\left\{
\begin{array}{rl}
>0,&j\in\mc{N}_i^{{\scriptsize \mbox{in}}},\\
0,&\mbox{otw.},
\end{array}
\right.
\qquad
\sum_{j=1}^na_{ij}=1,\forall i.
\end{equation*}
\begin{figure}[!h]
	\begin{center}
		\noindent
		\includegraphics[width=3in]{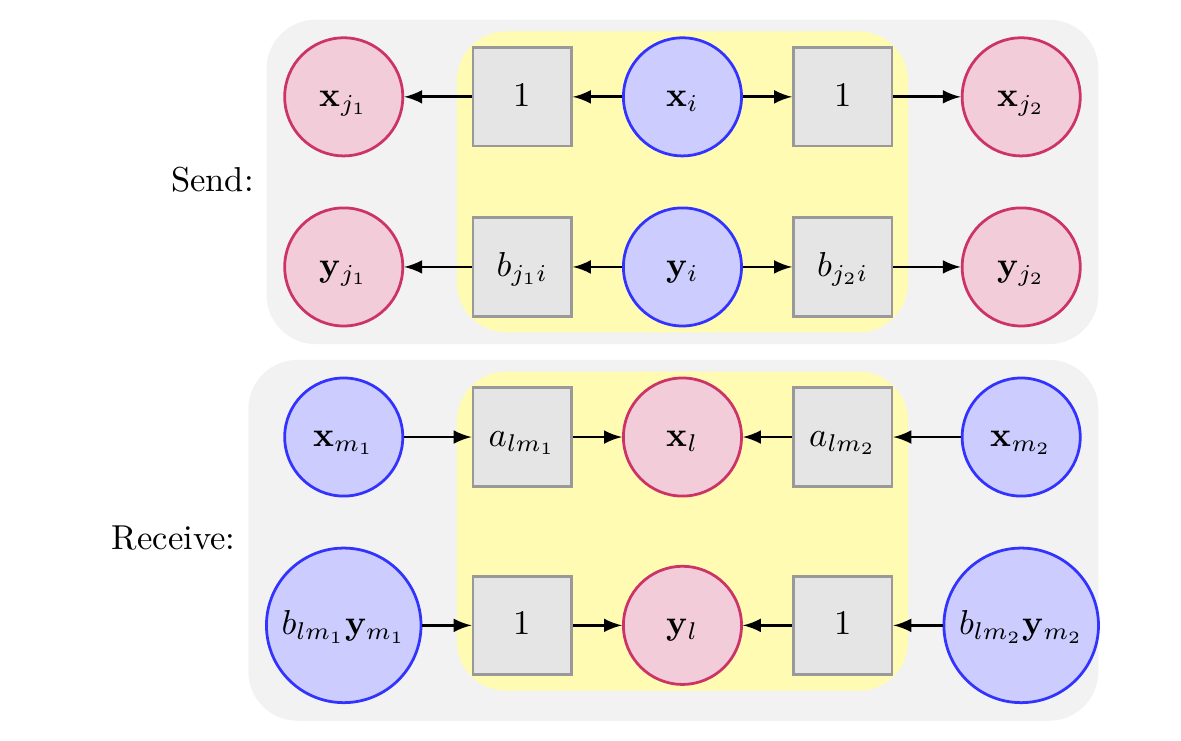}
		\caption{Illustration of the message passing between agents by Eq.~\eqref{alg1}.}\label{f_alg}
	\end{center}
\end{figure}
The diminishing step-size,~$\alpha_k\geq0$, satisfies the persistence conditions,~\cite{kushner2003stochastic,cc_lobel}:~$\sum_{k=0}^\infty\alpha_k=\infty$,~$\sum_{k=0}^\infty\alpha_k^2<\infty$. The scalar,~$\epsilon$, is a small positive number, which plays a key role in the convergence of the algorithm\footnote{Note that in the implementation of Eq.~\eqref{alg1}, each agent needs the knowledge of its out-neighbors. In a more restricted setting, e.g., a broadcast application where it may not be possible to know the out-neighbors, we may use~$b_{ij}=|\mc{N}_j^{{\scriptsize \mbox{out}}}|^{-1}$; thus, the implementation only requires knowing the out-degrees, see, e.g.,~\cite{opdirect_Nedic,opdirect_Nedic2} for similar assumptions.}. For an illustration of the message passing between agents in the implementation of Eq.~\eqref{alg1}, see Fig.~\ref{f_alg} on how agent~$i$ sends information to its out-neighbors and agent~$l$ receives information from its in-neighbors. In Fig.~\ref{f_alg}, the weights~$b_{j_1i}$ and~$b_{j_2i}$ are designed by agent~$i$, and satisfy~$b_{ii}+b_{j_1i}+b_{j_2i}=1$.  To analyze the algorithm, we denote~$\mb{z}_i^k\in\mbb{R}^p$,~$\mb{g}_i^k\in\mbb{R}^p$, and~$M\in\mbb{R}^{2n\times 2n}$ as follows:
\begin{align}
\mb{z}_i^k &= \left\{
\begin{array}{l r}
\mb{x}_i^k, & i\in\{1,...,n\},\\
\mb{y}_{i-n}^k, & ~~~~i\in\{n+1,...,2n\},
\end{array}
\right.\notag\\
\mb{g}_i^k &=
\left\{
\begin{array}{l r}
\nabla f_i(\mb{x}_i^k), & i\in\{1,...,n\},\\
0_p, & i\in\{n+1,...,2n\},\\
\end{array} \right.\notag\\
M&=\left[
\begin{array}{cc}
A & \epsilon I \\
I-A & B-\epsilon I \\
\end{array}
\right],\label{M}
\end{align}
where~$A=\{a_{ij}\}$ is row-stochastic,~$B=\{b_{ij}\}$ is column-stochastic.
Consequently, Eq.~\eqref{alg1} can be represented compactly as follows: for any~$i\in\{1,...,2n\}$, at~$k+1$th iteration,
\begin{align}\label{alg2}
\mb{z}_i^{k+1}=\sum_{j=1}^{2n}[M]_{ij}\mb{z}_j^{k}-\alpha_{k}\mb{g}_i^k.
\end{align}
We refer to the iterative relation in Eq.~\eqref{alg2} as the Directed-Distributed Gradient Descent (D-DGD) method, since it has the same form as DGD except the dimension doubles due to a new weight matrix~$M\in\mbb{R}^{2n\times 2n}$ as defined in Eq.~\eqref{M}. It is worth mentioning that even though Eq.~\eqref{alg2} looks similar to DGD,~\cite{uc_Nedic}, the convergence analysis of D-DGD does not exactly follow that of DGD. This is because the weight matrix,~$M$, has negative entries. Besides,~$M$ is not a doubly-stochastic matrix, i.e., the row sum is not~$1$. Hence, the tools in the analysis of DGD are not applicable, e.g.,~$\|\sum_{j}[M]_{ij}\mb{z}_j-\mb{x}^*\|\leq\sum_{j}[M]_{ij}\|\mb{z}_j-\mb{x}^*\|$ does not necessarily hold because~$[M]_{ij}$ are not non-negative. In next section, we prove the convergence of D-DGD.

\section{Convergence Analysis}\label{s3}
The convergence analysis of D-DGD can be divided into two parts. In the first part, we discuss the \emph{consensus property} of D-DGD, i.e., we capture the decrease in~$\left\|\mb{z}_i^k-\overline{\mb{z}}^k\right\|$ for~$i\in\{1,...,n\}$, as the D-DGD progresses, where we define~$\overline{\mb{z}}^k$ as the accumulation point:
\begin{align}\label{ac_point}
\overline{\mb{z}}^k\triangleq\frac{1}{n}\sum_{j=1}^{2n}\mb{z}_i^k=\frac{1}{n}\sum_{j=1}^{n}\mb{x}_i^k+\frac{1}{n}\sum_{j=1}^{n}\mb{y}_i^k.
\end{align}
The decrease in~$\left\|\mb{z}_i^k-\overline{\mb{z}}^k\right\|$ reveals that all agents approach a common accumulation point. We then show the \emph{optimality property} in the second part, i.e., the decrease in the difference between the function evaluated at the accumulation point and the optimal solution,~$f(\overline{\mb{z}}^k)-f(\mb{x}^*)$. We combine the two parts to establish the convergence.

\subsection{Consensus Property}
To show the consensus property, we study the convergence behavior of the weight matrices,~$M^k$, in Eq.~\eqref{M} as~$k$ goes to infinity. We use an existing results on such matrices~$M$, based on which we show the convergence behavior as well as the convergence rate. We borrow the following from~\cite{ac_Cai1}.
\begin{lem}\label{lem_M1}
(Cai~\textit{et al}.~\cite{ac_Cai1})  Assume the graph is strongly-connected.~$M$ is the weighting matrix defined in Eq.~\eqref{M}, and the constant~$\epsilon$ in~$M$ satisfies~$\epsilon\in(0,\Upsilon)$, where~$\Upsilon:=\frac{1}{(20+8n)^n}(1-|\lambda_3|)^n$, where~$\lambda_3$ is the third largest eigenvalue of~$M$ in Eq.~\eqref{M} by setting~$\epsilon=0$. Then the weighting matrix,~$M$, defined in Eq.~\eqref{M}, has a simple eigenvalue~$1$ and all other eigenvalues have magnitude smaller than one.
\end{lem}
\noindent Based on Lemma~\ref{lem_M1}, we now provide the convergence behavior as well as the convergence rate of the weight matrix,~$M$.
\begin{lem}\label{lem_M2}
	Assume that the network is strongly-connected, and~$M$ is the weight matrix that defined in Eq.~\eqref{M}.Then,
	\begin{enumerate}[label=(\alph*)]
		\item The sequence of~$\left\{M^k\right\}$, as~$k$ goes to infinity, converges to the following limit:
		\begin{align}
		\lim_{k\rightarrow\infty}M^k=\left[
		\begin{array}{cc}
		\frac{\mb{1}_n\mb{1}_n^\top}{n} & \frac{\mb{1}_n\mb{1}_n^\top}{n} \\
		0 & 0 \\
		\end{array}
		\right];\nonumber
		\end{align}
		\item For all~$i,j\in\mc{V}$, the entries~$\left[M^k\right]_{ij}$ converge to their limits as~$k\rightarrow\infty$ at a geometric rate, i.e., there exist bounded constants,~$\Gamma\in\mbb{R}$, and~$0<\gamma<1$, such that
		\begin{align}
		\left\|M^k-\left[
		\begin{array}{cc}
		\frac{\mb{1}_n\mb{1}_n^\top}{n} & \frac{\mb{1}_n\mb{1}_n^\top}{n} \\
		0 & 0 \\
		\end{array}
		\right]\right\|_{\infty}\leq\Gamma\gamma^k.\nonumber
		\end{align}
	\end{enumerate}
\end{lem}
\begin{pf}
Note that the sum of each column of~$M$ equals one, so~$1$ is an eigenvalue of~$M$ with a corresponding left (row) eigenvector~$[\mb{1}_n^\top~ \mb{1}_n^\top]$. We further have~$M[\mb{1}_n^\top~\mb{0}_n^\top]^\top=[\mb{1}_n^\top~\mb{0}_n^\top]^\top$, so~$[\mb{1}_n^\top~\mb{0}_n^\top]^\top$ is a right (column) eigenvector corresponding to the eigenvalue~$1$. According to Lemma~\ref{lem_M1},~$1$ is a simple eigenvalue of~$M$ and all other eigenvalues have magnitude smaller than one. We represent~$M^k$ in the Jordan canonical form for some~$P_i$ and~$Q_i$
\begin{align}\label{e1_lem_M2}
M^k=\frac{1}{n}[\mb{1}_n^\top~ \mb{0}_n^\top]^\top[\mb{1}_n^\top~\mb{1}_n^\top]+\sum_{i=2}^nP_iJ_i^kQ_i,
\end{align}
where the diagonal entries in~$J_i$ are smaller than one in magnitude for all~$i$. The statement (a) follows by noting that~$\lim_{k\rightarrow\infty}J_i^k=0$, for all~$i$.\\

From  Eq.~\eqref{e1_lem_M2}, and with the fact that all eigenvalues of~$M$ except~$1$ have magnitude smaller than one, there exist some bounded constants,~$\Gamma$ and~$\gamma\in(0,1)$, such that
\begin{align}
\left\|M^k-\left[
\begin{array}{cc}
\frac{\mb{1}_n\mb{1}_n^\top}{n} & \frac{\mb{1}_n\mb{1}_n^\top}{n} \\
0 & 0 \\
\end{array}
\right]\right\|&=\left\|\sum_{i=2}^nP_iJ_i^kQ_i\right\|,\nonumber\\
&\leq\sum_{i=2}^n\left\|P_i\right\|\left\|Q_i\right\|\left\|J_i^k\right\|\leq\Gamma\gamma^k\nonumber,
\end{align}
from which we get the desired result. \QEDB
\end{pf}

Using the result from Lemma~\ref{lem_M1}, Lemma~\ref{lem_M2} shows the convergence behavior of the power of the weight matrix, and further show that its convergence is bounded by a geometric rate. Lemma~\ref{lem_M2} plays a key role in proving the consensus properties of D-DGD. Based on Lemma~\ref{lem_M2}, we bound the difference between agent estimates in the following lemma. More specifically, we show that the agent estimates,~$\mb{x}_i^k$, approaches the accumulation point,~$\overline{\mb{z}}^k$, and the auxiliary variable,~$\mb{y}_i^k$, goes to~$\mb{0}_n$, where~$\overline{\mb{z}}^k$ is defined in Eq.~\eqref{ac_point}.
\begin{lem}\label{lem_consensus}
Let the Assumptions A1 hold. Let~$\left\{\mb{z}_i^k\right\}$ be the sequence over~$k$ generated by the D-DGD algorithm, Eq.~\eqref{alg2}. Then, there exist some bounded constants,~$\Gamma$ and~$0<\gamma<1$, such that:
\begin{enumerate}[label=(\alph*)]
\item for~$1\leq i\leq n$, and~$k\geq1$,
		\begin{align}
		\left\|\mb{z}_i^k-\overline{\mb{z}}^k\right\|\leq&\Gamma\gamma^k\sum_{j=1}^{2n}\left\|\mb{z}_j^0\right\|+n\Gamma D\sum_{r=1}^{k-1}\gamma^{k-r}\alpha_{r-1}\nonumber\\
&+2D\alpha_{k-1};\nonumber
		\end{align}

\item for~$n+1\leq i\leq 2n$, and~$k\geq1$,
		\begin{align}
		\left\|\mb{z}_i^k\right\|\leq&\Gamma\gamma^k\sum_{j=1}^{2n}\left\|\mb{z}_j^0\right\|+n\Gamma D\sum_{r=1}^{k-1}\gamma^{k-r}\alpha_{r-1}.\nonumber
		\end{align}
\end{enumerate}
\end{lem}
\begin{pf}
	For any~$k\geq 1$, we write Eq.~\eqref{alg2} recursively
	\begin{align}\label{eq1_lem_consensus}
	\mb{z}_i^{k}=&\sum_{j=1}^{2n}[M^k]_{ij}\mb{z}_j^0-\sum_{r=1}^{k-1}\sum_{j=1}^{2n}[M^{k-r}]_{ij}\alpha_{r-1}\mb{g}_j^{r-1}\nonumber\\
&-\alpha_{k-1}\mb{g}_i^{k-1}.
	\end{align}	
	Since every column of~$M$ sums up to one, we have for any~$r$ $\sum_{i=1}^{2n}[M^{r}]_{ij}=1$. Considering the recursive relation of~$\mb{z}_i^{k}$ in Eq.~\eqref{eq1_lem_consensus}, we obtain that~$\overline{\mb{z}}^k$ can be represented as
	\begin{align}\label{eq2_lem_consensus}
	\overline{\mb{z}}^k&=\sum_{j=1}^{2n}\frac{1}{n}\mb{z}_j^0-\sum_{r=1}^{k-1}\sum_{j=1}^{2n}\frac{1}{n}\alpha_{r-1}\mb{g}_j^{r-1}-\frac{1}{n}\sum_{j=1}^{2n}\alpha_{k-1}\mb{g}_j^{k-1}.
	\end{align}	
	Subtracting Eq.~\eqref{eq2_lem_consensus} from~\eqref{eq1_lem_consensus} and taking the norm, we obtain that for~$1\leq i\leq n$,
	\begin{align}\label{eq3_lem_consensus}
	&\left\|\mb{z}_i^{k}-\overline{\mb{z}}^k\right\|\leq\sum_{j=1}^{2n}\left\|[M^k]_{ij}-\frac{1}{n}\right\|\left\|\mb{z}_j^0\right\|\nonumber\\
	&+\sum_{r=1}^{k-1}\sum_{j=1}^{n}\left\|[M^{k-r}]_{ij}-\frac{1}{n}\right\|\alpha_{r-1}\left\|\nabla f_j(\mb{x}_j^{r-1})\right\|\nonumber\\
	&+\alpha_{k-1}\left\|\nabla f_i(\mb{x}_i^{k-1})\right\|+\frac{1}{n}\sum_{j=1}^n\alpha_{k-1}\left\|\nabla f_j(\mb{x}_j^{k-1})\right\|.
	\end{align}
	The proof of part (a) follows by applying the result of Lemma~\ref{lem_M2} to Eq.~\eqref{eq3_lem_consensus} and noticing that the gradient is bounded by a constant~$D$. Similarly, by taking the norm of Eq.~\eqref{eq1_lem_consensus}, we obtain that for~$n+1\leq i\leq 2n$,
	\begin{align}
	\left\|\mb{z}_i^{k}\right\|&\leq\sum_{j=1}^{2n}\left\|[M^k]_{ij}\right\|\left\|\mb{z}_j^0\right\|\nonumber\\
&+\sum_{r=1}^{k-1}\sum_{j=1}^{n}\left\|[M^{k-r}]_{ij}\right\|\alpha_{r-1}\left\|\nabla f_j(\mb{x}_j^{r-1})\right\|.\nonumber
	\end{align}
	The proof of part (b) follows by applying the result of Lemma \ref{lem_M2} to the preceding relation and considering the boundedness of gradient in Assumption~\ref{asp}(e).\QEDB
\end{pf}

Using the above lemma, we now draw our first conclusion on the consensus property at the agents. Proposition~\ref{prop_consensus} reveals that all agents asymptotically reach consensus.
\begin{prop}\label{prop_consensus}
Let the Assumptions A1 hold. Let~$\left\{\mb{z}_i^k\right\}$ be the sequence over~$k$ generated by the D-DGD algorithm, Eq.~\eqref{alg2}. Then,~$\mb{z}_i^k$ satisfies
\begin{enumerate}[label=(\alph*)]
\item for~$1\leq i\leq n$,
		\begin{align}
			\sum_{k=1}^{\infty}\alpha_k\left\|\mb{z}_i^k-\overline{\mb{z}}^k\right\|<\infty;\nonumber
		\end{align}

\item for~$n+1\leq i\leq 2n$,
		\begin{align}
			\sum_{k=1}^{\infty}\alpha_k\left\|\mb{z}_i^k\right\|<\infty.\nonumber
		\end{align}
\end{enumerate}
\end{prop}
\begin{pf}
	Based on the result of Lemma~\ref{lem_consensus}(a), we obtain, for~$1\leq i\leq n$,
	\begin{align}\label{prop_Consensus_eqa}
	\sum_{k=1}^{K}&\alpha_k\left\|\mb{z}_i^k-\overline{\mb{z}}^k\right\|\leq \Gamma\left(\sum_{j=1}^{2n}\left\|\mb{z}_j^0\right\|\right)\sum_{k=1}^K\alpha_k\gamma^{k}\nonumber\\
	&+n\Gamma D\sum_{k=1}^K\sum_{r=1}^{k-1}\gamma^{(k-r)}\alpha_k\alpha_{r-1}+2D\sum_{k=0}^{K-1}\alpha_k^2.
	\end{align}
	With the basic inequality~$ab\leq\frac{1}{2}(a^2+b^2)$,~$a,b\in\mbb{R}$, we have:
	\begin{align}\nonumber
	2\sum_{k=1}^K\alpha_k\gamma^{k}\leq\sum_{k=1}^K\left[\alpha_k^2+\gamma^{2k}\right]\leq\sum_{k=1}^K\alpha_k^2+\frac{1}{1-\gamma^2};
	\end{align}
and
	\begin{align}\nonumber
	\sum_{k=1}^K&\sum_{r=1}^{k-1}\gamma^{(k-r)}\alpha_k\alpha_{r-1}\leq\frac{1}{2}\sum_{k=1}^K\alpha_k^2\sum_{r=1}^{k-1}\gamma^{(k-r)}\\
	&+\frac{1}{2}\sum_{r=1}^{K-1}(\alpha_{r-1})^2\sum_{k=r+1}^{K}\gamma^{(k-r)}\leq\frac{1}{1-\gamma}\sum_{k=1}^K\alpha_k^2.\nonumber
	\end{align}
	The proof of part (a) follows by applying the preceding relations to Eq.~\eqref{prop_Consensus_eqa} along with~$\sum_{k=0}^K\alpha_k^2<\infty$ as~$K\rightarrow\infty$.
    Following the same spirit in the proof of part (b), we can reach the conclusion of part (b).\QEDB
\end{pf}
Since~$\sum_{k=1}^\infty\alpha_k=\infty$, Proposition~\ref{prop_consensus} shows that all agents reach consensus at the accumulation point,~$\overline{\mb{z}}^k$, asymptotically, i.e.,~for all~$1\leq i\leq n$,~$1\leq j\leq n$,
\begin{align}\label{eq_lim1}
\lim_{k\rightarrow\infty}\mb{z}_i^k=\lim_{k\rightarrow\infty}\overline{\mb{z}}^k=\lim_{k\rightarrow\infty}\mb{z}_j^k,
\end{align}
and for~$n+1\leq i\leq 2n$, the states,~$\mb{z}_i^k$, asymptotically, converge to zero, i.e., for~$n+1\leq i\leq 2n$,
\begin{align}\label{eq_lim2}
\lim_{k\rightarrow\infty}\mb{z}_i^k=0.
\end{align}
We next show how the accumulation point,~$\overline{\mb{z}}^k$, approaches the optima,~$\mb{x}^*$, as D-DGD progresses.

\subsection{Optimality Property}
The following lemma gives an upper bound on the difference between the objective evaluated at the accumulation point,~$f(\overline{\mb{z}}^k)$, and the optimal objective value,~$f^*$.
\begin{lem}\label{lem_opt1}
Let the Assumptions A1 hold. Let~$\left\{\mb{z}_i^k\right\}$ be the sequence over~$k$ generated by the D-DGD algorithm, Eq.~\eqref{alg2}. Then,
\begin{align}\label{lem_opt1_eq}
	2\sum_{k=0}^\infty\alpha_k\left(f(\overline{\mb{z}}^k)-f^*\right)&\leq n\left\|\overline{\mb{z}}^0-\mb{x}^*\right\|^2+nD^2\sum_{k=0}^\infty\alpha_k^2\nonumber\\
	&+\frac{4D}{n}\sum_{i=1}^{n}\sum_{k=0}^\infty\alpha_k\left\|\mb{z}_i^k-\overline{\mb{z}}^k\right\|.
	\end{align}	
	\end{lem}
\begin{pf}
	Consider Eq.~\eqref{alg2} and the fact that each column of~$M$ sums to one, we have
	\begin{align}
	\overline{\mb{z}}^{k+1}&=\frac{1}{n}\sum_{j=1}^{2n}\left[\sum_{i=1}^{2n}[M]_{ij}\right]\mb{z}_j^{k}-\alpha_{k}\frac{1}{n}\sum_{i=1}^{2n}\mb{g}_i^k,\nonumber\\
	&=\overline{\mb{z}}^k-\frac{\alpha_k}{n}\sum_{i=1}^{n}\nabla f_i(\mb{z}_i^k).\nonumber
	\end{align}
	Therefore, we obtain that
	\begin{align}
	\left\|\overline{\mb{z}}^{k+1}-\mb{x}^*\right\|^2&=\left\|\overline{\mb{z}}^k-\mb{x}^*\right\|^2+\left\|\frac{\alpha_k}{n}\sum_{i=1}^{n}\nabla f_i(\mb{z}_i^k)\right\|^2\nonumber\\
	&-2\frac{\alpha_k}{n}\sum_{i=1}^{n}\left\langle\overline{\mb{z}}^k-\mb{x}^*,\nabla f_i(\mb{z}_i^k)\right\rangle\label{mr_eq1}.
	\end{align}
	Denote $\nabla f_i^k=\nabla f_i(\mb{z}_i^k)$. Since $\|\nabla f_i^k\|\leq D$, we have 
	\begin{align}
	&\left\langle\overline{\mb{z}}^{k}-\mb{x}^{*},\nabla f_i^k\right\rangle=\left\langle\overline{\mb{z}}^{k}-\mb{z}_i^{k},\nabla f_i^k\right\rangle+\left\langle\mb{z}_i^{k}-\mb{x}^{*},\nabla f_i^k\right\rangle\nonumber\\
	&\geq \left\langle\overline{\mb{z}}^{k}-\mb{z}_i^{k},\nabla f_i^k\right\rangle+f_i(\mb{z}_i^k)-f_i(\mb{x}^*)\nonumber\\
	&\geq-D\left\|\mb{z}_i^{k}-\overline{\mb{z}}^{k}\right\|+f_i(\mb{z}_i^k)-f_i(\overline{\mb{z}}^{k})+f_i(\overline{\mb{z}}^{k})-f_i(\mb{x}^*)\nonumber\\
	&\geq-2D\left\|\mb{z}_i^{k}-\overline{\mb{z}}^{k}\right\|+f_i(\overline{\mb{z}}^{k})-f_i(\mb{x}^*).\label{mr_eq2}
	\end{align}	
	By substituting Eq.~\eqref{mr_eq2} in Eq.~\eqref{mr_eq1}, and rearranging the terms, we obtain that
	\begin{align}
	2\alpha_k\left(f(\overline{\mb{z}}^k)-f^*\right)&\leq n\left\|\overline{\mb{z}}^k-\mb{x}^*\right\|^2-n\left\|\overline{\mb{z}}^{k+1}-\mb{x}^*\right\|^2\nonumber\\
	&+nD^2\alpha_k^2+\frac{4D}{n}\sum_{i=1}^{n}\alpha_k\left\|\mb{z}_i^k-\overline{\mb{z}}^k\right\|.\label{mr_eq3}
	\end{align}
	The desired result is achieved by summing Eq.~\eqref{mr_eq3} over time from $k=0$ to $\infty$.	
\QEDB
\end{pf}
We are ready to present the main result of this paper, by combining all the preceding results.
\begin{thm}
	Let the Assumptions A1 hold. Let~$\left\{\mb{z}_i^k\right\}$ be the sequence over~$k$ generated by the D-DGD algorithm, Eq.~\eqref{alg2}. Then, for any agent $i$, we have 
	\begin{align}
	\lim_{k\rightarrow\infty}f(\mb{z}_i^k)=f^*.\nonumber
	\end{align}
\end{thm}
\begin{pf}
	Since that the step-size follows that $\sum_{k=0}^\infty\alpha_k^2<\infty$, and $\sum_{k=0}^\infty\alpha_k\|\mb{z}_i^{k}-\overline{\mb{z}}^{k}\|<\infty$ from Lemma \ref{prop_consensus}, we obtain from Eq.~\eqref{lem_opt1_eq} that
	\begin{align}
	2\sum_{k=0}^\infty\alpha_k\left(f(\overline{\mb{z}}^{k})-f^*\right)<\infty,
	\end{align}
	which reveals that $\lim_{k\rightarrow\infty}f(\overline{\mb{z}}^{k})=f^*$, by considering that $\sum_{k=0}^\infty\alpha_k=\infty$. In Eq.~\eqref{eq_lim1}, we have already shown that $\lim_{k\rightarrow\infty}\mb{z}_i^k=\lim_{k\rightarrow\infty}\overline{\mb{z}}^k$. Therefore, we obtain the desired result.\QEDB
\end{pf}
\section{Convergence Rate}\label{s4}
In this section, we show the convergence rate of D-DGD. Let~$f_m:=\min_{k}f(\overline{\mb{z}}^k)$, we have
\begin{align}\label{rate_ineq}
(f_m-f^*)\sum_{k=0}^K\alpha_k\leq\sum_{k=0}^K\alpha_k(f(\overline{\mb{z}}^k)-f^*)
\end{align}

By combining Eqs.~\eqref{prop_Consensus_eqa},~\eqref{lem_opt1_eq} and~\eqref{rate_ineq}, it can be verified that Eq.~\eqref{lem_opt1_eq} can be represented in the following form:
\begin{align}
(f_m-f^*)\sum_{k=0}^K\alpha_k\leq C_1+C_2\sum_{k=0}^K\alpha_k^2,\nonumber
\end{align}
or equivalently,
\begin{align}\label{rate}
(f_m-f^*)\leq\frac{ C_1}{\sum_{k=0}^K\alpha_k}+\frac{C_2\sum_{k=0}^K\alpha_k^2}{\sum_{k=0}^K\alpha_k},
\end{align}
where the constants,~$C_1$ and~$C_2$, are given by
\begin{align}
C_1=&\frac{n}{2}\left\|\overline{\mb{z}}^0-\mb{x}^*\right\|^2-\frac{n}{2}\left\|\overline{\mb{z}}^{K+1}-\mb{x}^*\right\|^2\nonumber\\
&+D\Gamma\sum_{j=1}^{2n}\left\|\mb{z}_j^0\right\|\frac{1}{1-\gamma^2},\nonumber\\
C_2=&\frac{nD^2}{2}+4D^2+D\Gamma\sum_{j=1}^{2n}\left\|\mb{z}_j^0\right\|+\frac{2D^2\Gamma}{1-\gamma}.\nonumber
\end{align}
Eq.~\eqref{rate} actually has the same form as the \mbox{equations} in analyzing the convergence rate of DGD (recall, e.g.,~\cite{uc_Nedic}). In particular, when~$\alpha_k=k^{-1/2}$, the first term in Eq.~\eqref{rate}  leads to
\begin{align}
\frac{ C_1}{\sum_{k=0}^K\alpha_k}=C_1\frac{1/2}{K^{1/2}-1}=O\left(\frac{1}{\sqrt{K}}\right),\nonumber
\end{align}
while the second term in Eq.~\eqref{rate} leads to
\begin{align}
\frac{C_2\sum_{k=0}^K\alpha_k^2}{\sum_{k=0}^K\alpha_k}=C_2\frac{
	\ln K}{2(\sqrt{K}-1)}=O\left(\frac{\ln K}{\sqrt{K}}\right).\nonumber
\end{align}
It can be observed that the second term dominates, and the overall convergence rate is~$O\left(\frac{\ln k}{\sqrt{k}}\right)$. As a result, D-DGD has the same convergence rate as DGD. The restriction of directed graph does not effect the speed.
\section{Numerical Experiment}\label{s5}
\begin{figure}[!htb]
\begin{minipage}[b]{.32\linewidth}
  \centering
  \centerline{\includegraphics[width=2.8cm]{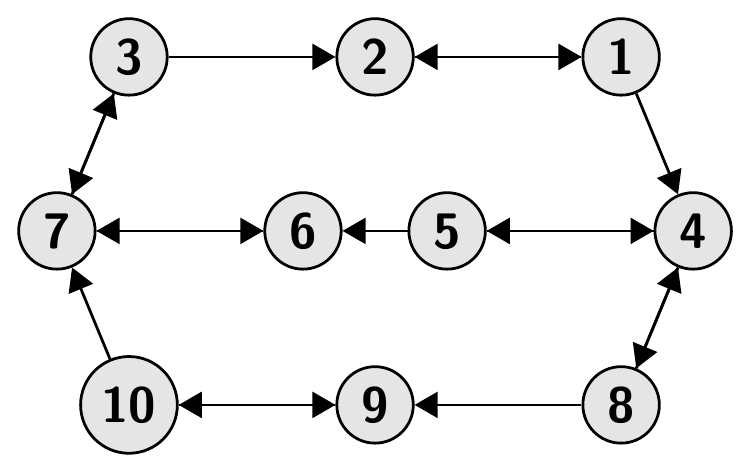}}
  \centerline{$\mc{G}_a$}\medskip
\end{minipage}
\hfill
\begin{minipage}[b]{.32\linewidth}
  \centering
  \centerline{\includegraphics[width=2.8cm]{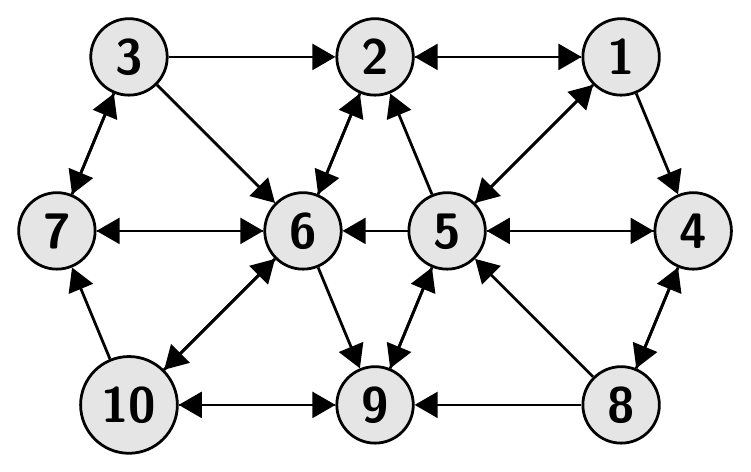}}
  \centerline{$\mc{G}_b$}\medskip
\end{minipage}
\hfill
\begin{minipage}[b]{.32\linewidth}
  \centering
  \centerline{\includegraphics[width=2.8cm]{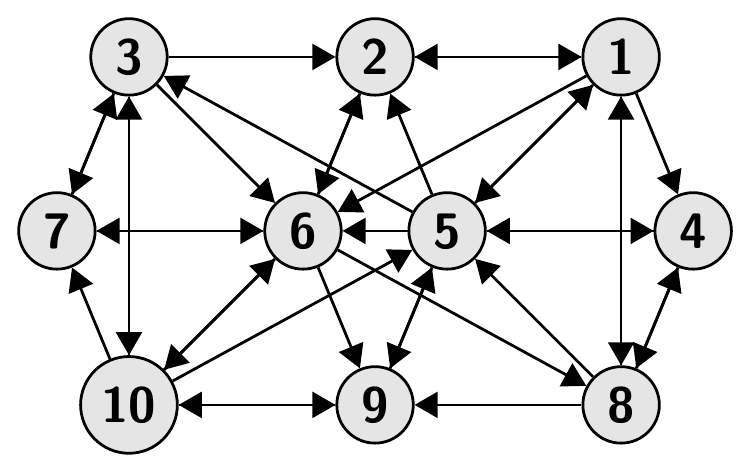}}
  \centerline{$\mc{G}_c$}\medskip
\end{minipage}
\caption{Three examples of strongly-connected but non-balanced digraphs.}
\label{graph}
\end{figure}

We consider a distributed least squares problem in a directed graph: each agent owns a private objective function,~$\mb{s}_i=R_i\mb{x}+\mb{n}_i$, where~$\mb{s}_i\in\mbb{R}^{m_i}$ and~$R_i\in\mbb{R}^{m_i\times p}$ are measured data,~$\mb{x}\in\mbb{R}^p$ is unknown states, and~$\mb{n}_i\in\mbb{R}^{m_i}$ is random unknown noise. The goal is to estimate~$\mb{x}$. This problem can be formulated as a distributed optimization problem solving
\begin{align}
\mbox{min }f(\mb{x})=\frac{1}{n}\sum_{i=1}^n\left\|R_i\mb{x}-\mb{s}_i\right\|.\nonumber
\end{align}
We consider the network topology as the digraphs shown in Fig.~\ref{graph}. We employ identical setting and graphs as~\cite{ac_Cai1}. In~\cite{ac_Cai1}, the value of ~$\epsilon=0.7$ is chosen for each~$\mc{G}_a,\mc{G}_b,\mc{G}_c$.

\begin{figure}[!h]
	\begin{center}
		\noindent
		\includegraphics[width=3.5in]{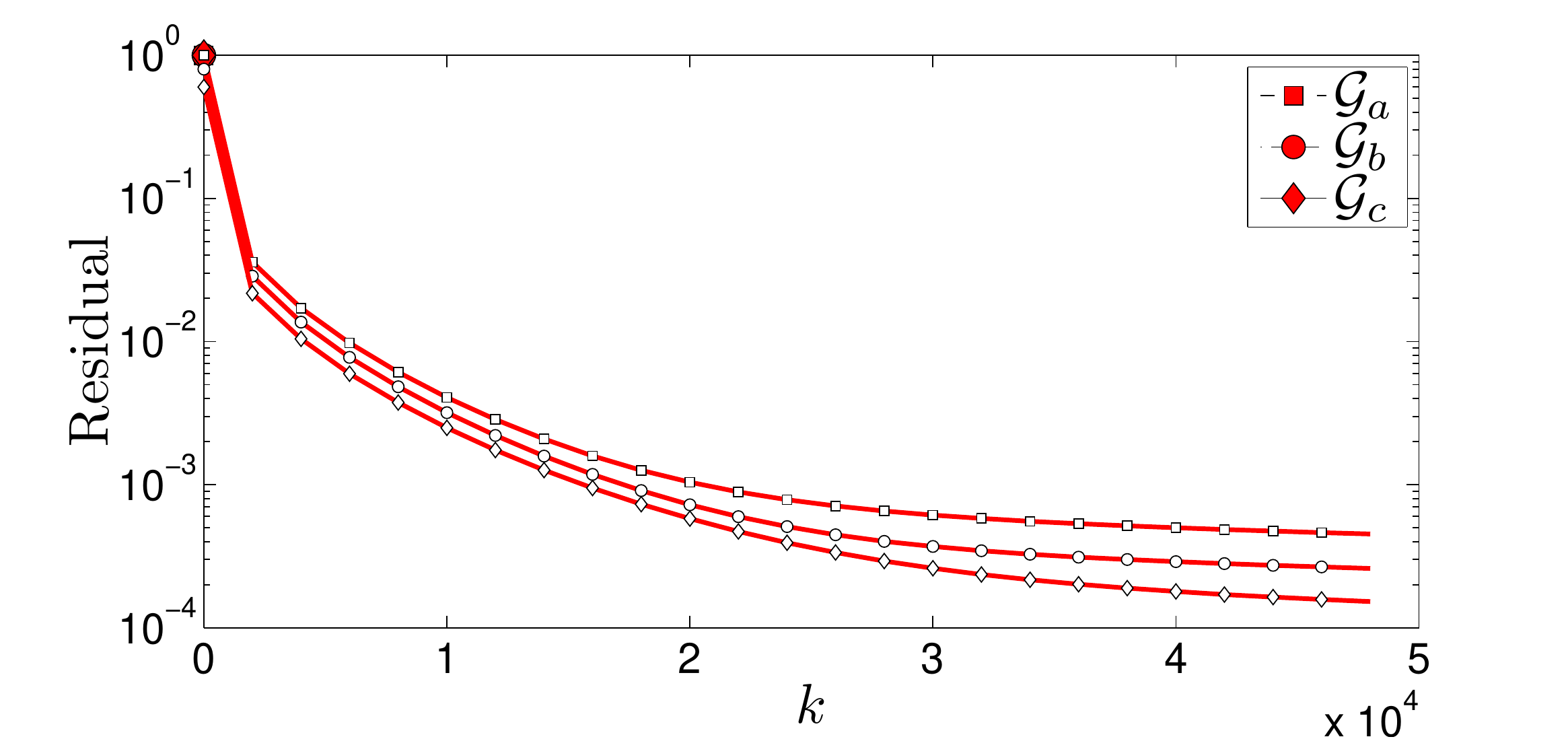}
		\caption{Plot of residuals~$\frac{\left\|\mb{x}_k-\mb{x}^*\right\|_F}{\left\|\mb{x}_0-\mb{x}^*\right\|_F}$ for digraph~$\mc{G}_a,\mc{G}_b,\mc{G}_c$ as D-DGD progresses.}\label{diff_graph}
	\end{center}
\end{figure}

Fig.~\ref{diff_graph} shows the convergence of the D-DGD algorithm for three digraphs displayed in Fig.~\ref{graph}. Once the weight matrix,~$M$, defined in Eq.~\eqref{M}, converges, the D-DGD ensures the convergence. Moreover, it can be observed that the residuals decrease faster as the number of edges increases, from~$\mc{G}_a$ to~$\mc{G}_c$. This indicates faster convergence when there are more communication channels available for information exchange.

\begin{figure}[!h]
	\begin{center}
		\noindent
		\includegraphics[width=3.5in]{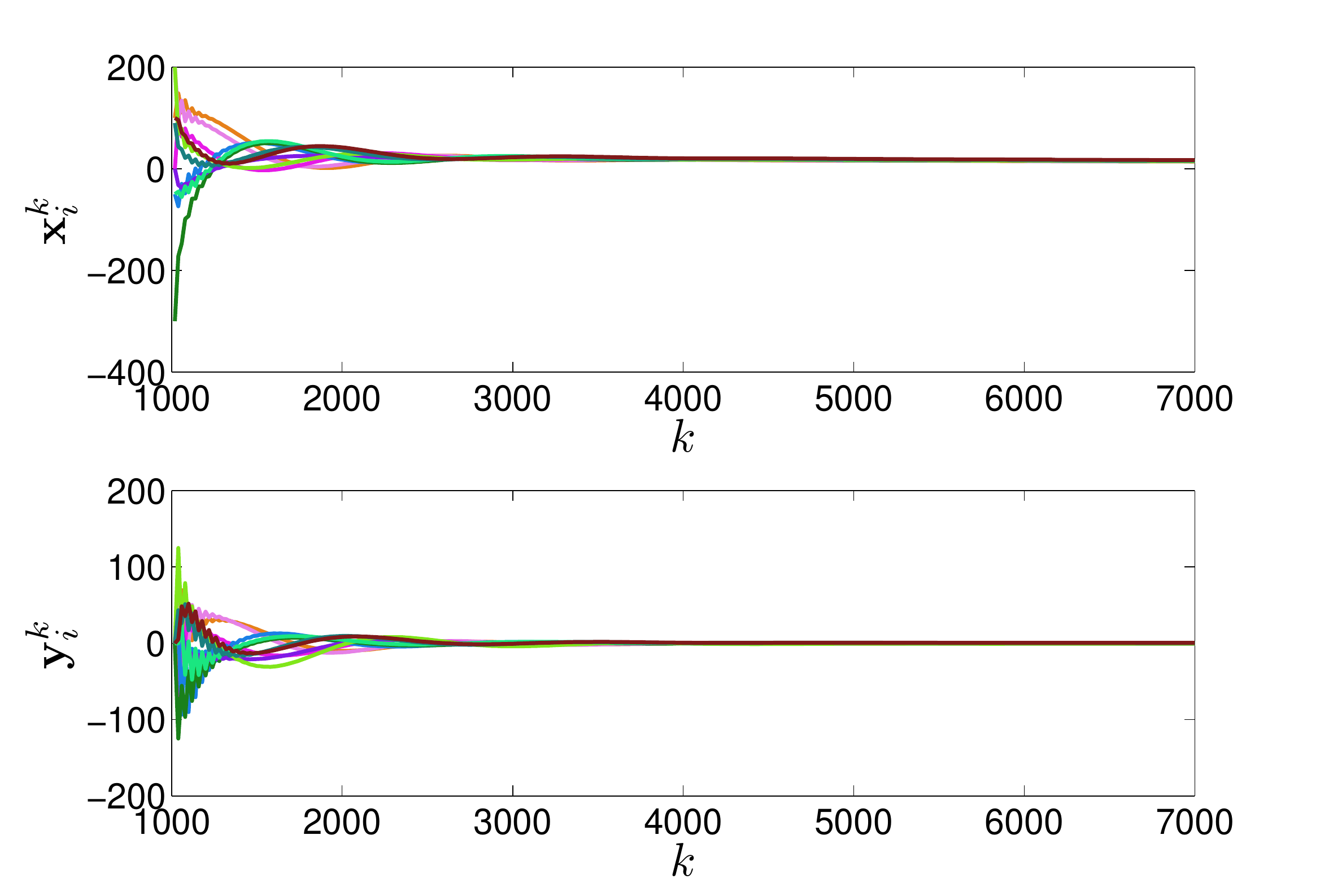}
		\caption{Sample paths of states,~$\mb{x}_i^k$, and~$\mb{y}_i^k$, for all agents on digraphs~$\mc{G}_a$ with~$\epsilon=0.7$ as D-DGD progresses.}\label{path}
	\end{center}
\end{figure}

In Fig.~\ref{path}, we display the trajectories of both states,~$\mb{x}$ and~$\mb{y}$, when the D-DGD, Eq.~\eqref{alg2}, is applied on digraph~$\mc{G}_a$ with parameter~$\epsilon=0.7$. Recall that in Eqs.~\eqref{eq_lim1} and~\eqref{eq_lim2}, we have shown that as times,~$k$, goes to infinity, the state,~$\mb{x}_i^k$ of all agents will converges to a same accumulation point,~$\overline{\mb{z}}^k$, which is the optimal solution of the problem, and~$\mb{y}_i^k$ of all agents converges to zero, which are shown in Fig.~\ref{path}.

\begin{figure}[!h]
	\begin{center}
		\noindent
		\includegraphics[width=3.5in]{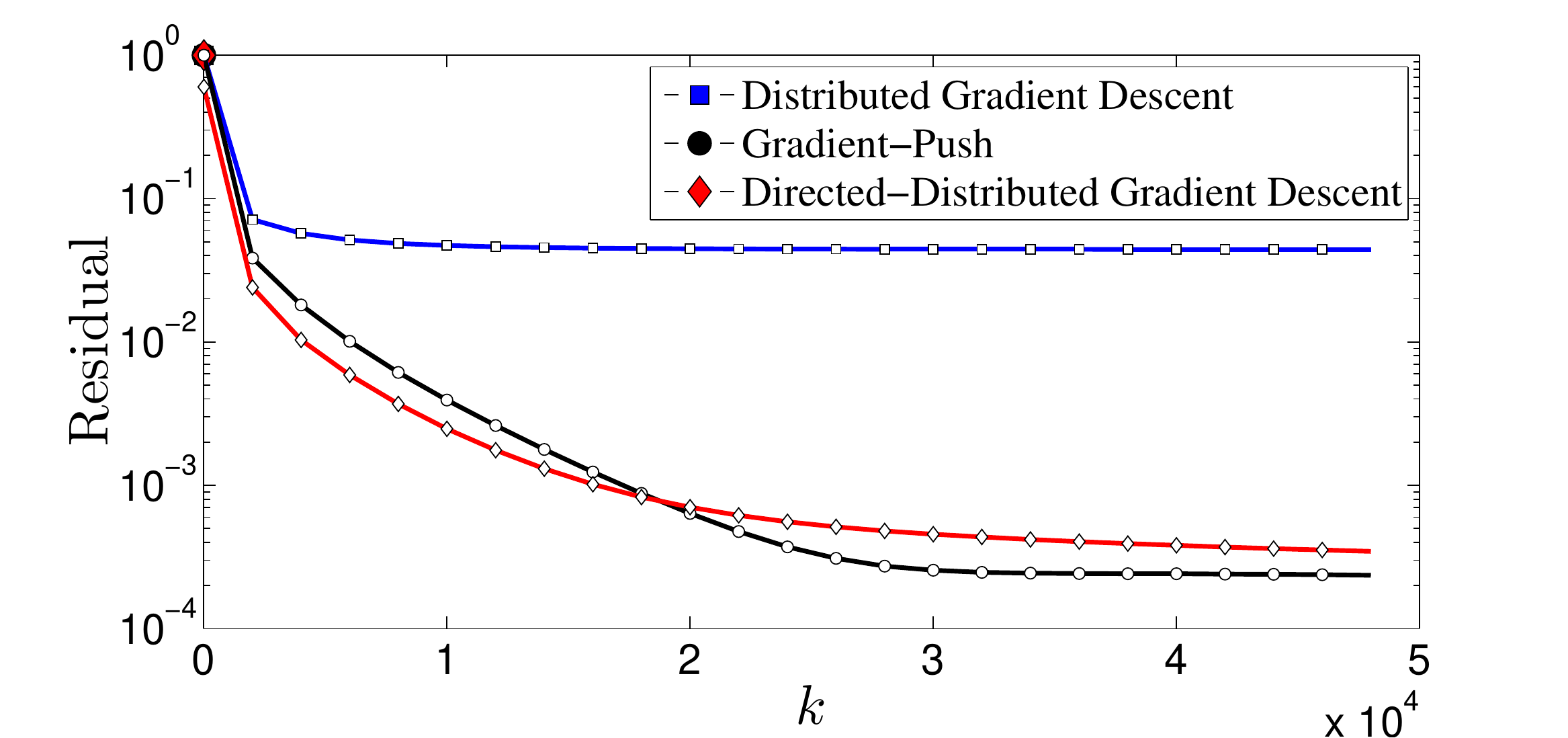}
		\caption{Plot of residuals~$\frac{\left\|\mb{x}_k-\mb{x}^*\right\|_F}{\left\|\mb{x}_0-\mb{x}^*\right\|_F}$ as (D-)DGD progresses.}\label{diff_alg}
	\end{center}
\end{figure}

In the next experiment, we compare the performance between the D-DGD and others distributed optimization algorithms over directed graphs. The red curve in Fig.~\ref{diff_alg} is the plot of residuals of D-DGD on~$\mc{G}_a$. In Fig.~\ref{diff_alg}, we also shown the convergence behavior of two other algorithms on the same digraph. The blue line is the plot of residuals with a DGD algorithm using a row-stochastic matrix.
As we have discussed is Section~\ref{s2}, when the weight matrix is restricted to be row-stochastic, DGD actually minimizes a new objective function~$\widehat{f}(\mb{x})=\sum_{i=1}^n\pi_if_i(\mb{x})$ where~$\boldsymbol{\pi}=\{\pi_i\}$ is the left eigenvector of the weight matrix corresponding to eigenvalue~$1$. So it does not converge to the true~$\mb{x}^*$. The black curve shows the convergence behavior of the gradient-push algorithm, proposed in~\cite{opdirect_Nedic,opdirect_Nedic2}. Our algorithm has the same convergence rate as the gradient-push algorithm, which is~$O\left(\frac{\ln k}{\sqrt{k}}\right)$.
\section{Conclusion}\label{s6}
In this paper, we describe a distributed algorithm, called Directed-Distributed Gradient Descent (D-DGD), to solve the problem of minimizing a sum of convex objective functions over a~\emph{directed} graph. Existing distributed algorithms,~e.g.,~Distributed Gradient Descent (DGD), deal with the same problem under the assumption of undirected networks. The primary reason behind assuming the undirected graphs is to obtain a doubly-stochastic weight matrix. The row-stochasticity of the weight matrix guarantees that all agents reach consensus, while the column-stochasticity ensures optimality, i.e., each agent¡¯s local gradient contributes equally to the global objective. In a directed graph, however, it may not be possible to construct a doubly-stochastic weight matrix in a distributed manner. In each iteration of D-DGD, we simultaneously constructs a row-stochastic matrix and a column-stochastic matrix instead of only a doubly-stochastic matrix. The convergence of the new weight matrix, depending on the row-stochastic and column-stochastic matrices, ensures agents to reach both consensus and optimality. The analysis shows that the D-DGD converges at a rate of~$O(\frac{\ln k}{\sqrt{k}})$, where~$k$ is the number of iterations.



\end{document}